\magnification=\magstep1
\input amstex
\documentstyle{amsppt}
\voffset=-3pc
\loadbold
\loadmsbm
\loadeufm
\UseAMSsymbols
\baselineskip=12pt
\parskip=6pt
\def\var{\varepsilon}
\def\bC{\Bbb C}
\def\bN{\Bbb N}
\def\bR{\Bbb R}

\def\bZ{\Bbb Z}
\def\cO{\Cal O}

\def\cA{\Cal A}
\def\cE{\Cal E}

\def\cF{\Cal F}
\def\cG{\Cal G}
\def\cJ{\Cal J}
\def\cB{\Cal B}
\def\cS{\Cal S}
\def\cC{\Cal C}
\def\cK{\Cal K}
\def\fW{\frak W}
\def\cI{\vartheta}
\def\fU{\frak U}
\def\fV{\frak V}
\def\cU{\Cal U}

\def\id{\text{id}}
\def\Im{\text{Im }}

\def\Ker{\text{Ker }}

\def\Hom{\text{Hom}}
\def\bHom{\bold{Hom}}

\def\plain{\text{plain}}
\def\be{\bold e}
\def\bg{\bold g}
\NoBlackBoxes
\topmatter
\title Analytic Cohomology Groups of Infinite Dimensional Complex Manifolds\endtitle
\address Department of Mathematics,
Purdue University,
West Lafayette IN\ 47907-1395\endaddress
\author
\bf L\'aszl\'o Lempert\endauthor\footnote""{Research supported by NSF
grants DMS0700281 and 1162070\hfill\break}
\leftheadtext{Analytic Cohomology Groups}
\rightheadtext{L\'aszl\'o Lempert}
\document
\TagsOnRight
\abstract Given a cohesive sheaf $\cS$ over a complex Banach manifold $M$, we endow the cohomology groups $H^q(M,\cS)$ of $M$ and $H^q(\fU,\cS)$ of open covers $\fU$ of $M$ with a locally convex topology.
Under certain assumptions we prove that the canonical map $H^q(\fU,\cS)\to H^q (M,\cS)$ is an isomorphism of topological vector spaces.\endabstract
\endtopmatter

\head 1.\ Introduction\endhead

Sheaf theory is an indispensable tool in the study of complex manifolds.
For finite dimensional manifolds the relevant theory is that of coherent sheaves, and for Banach manifolds with Patyi we developed in \cite{LP} an analogous theory of cohesive sheaves.
In \cite{LP} we primarily dealt with manifolds modeled on Banach spaces with an unconditional basis, but in light of Patyi's more recent work \cite{P2}, the results of \cite{LP} hold more generally, e.g.~for manifolds modeled on Banach spaces with a Schauder basis.
The focus of \cite{LP} was to prove that higher cohomology groups of cohesive sheaves $\cS\to M$ vanish, under a suitable convexity condition on the manifold $M$.
Our goal here is to introduce a topology on the space $\Gamma(M,\cS)$ of sections and on the cohomology groups $H^q(M,\cS)$ of rather general manifolds, and to prove a topological version of Leray's isomorphism theorem.
The topology on $H^q(M,\cS)$ is obtained as the direct limit of topologies on the \v Cech groups $H^q(\fU,\cS)$ corresponding to open covers $\fU$ of $M$, and our main result, Theorem 4.5, says that if the cohesive sheaf $\cS$ is separated in a certain sense, and elements of $\fU$ have a suitable convexity property ($\fU$ is a ``Stein'' cover), then the canonical map $H^q(\fU,\cS)\to H^q(M,\cS)$ is an isomorphism of topological vector spaces.

In the finite dimensional theory the importance of topology on cohomology groups has long been recognized, to study finiteness, compactly supported cohomology groups, analytic continuation, and embedding problems [Bu,CS,K,RR,RRV,S1-2].
By contrast, in infinite dimensions so far only $\Gamma(M,\cS)$ has been considered as a topological space, when $\cS$ is the sheaf of germs of holomorphic functions valued in some topological vector space $E$, so that $\Gamma(M,\cS)$ can be identified with the space $\cO^E(M)$ of holomorphic maps $M\to E$; see \cite{D2}.
When $\dim M<\infty$, there is only one reasonable topology on the groups $\Gamma(M,\cS)$ and $H^q(M,\cS)$.
Not so for infinite dimensional $M$.
There are several equally reasonable topologies on the spaces $\cO^E(M)$, and each gives rise to a topology on $H^q(M,\cS)$.
However, there is only one among these topologies for which a topological Leray isomorphism could be proved.
It is induced by the so--called $\tau_\delta$ or countable cover topology on $\cO^E(M)$, first studied by Coeur\'e and Nachbin in a somewhat lesser generality \cite{C,N}.
The advantage of this topology is that it is the direct limit of Banach space topologies.

Since the canonical map $H^q(\fU,\cS)\to H^q(M,\cS)$ is always continuous, and under a suitable convexity condition on the cover $\fU$ 
it is bijective, to prove our main result we need to show that it is open as well. This will follow from two theorems.
The first of these two is a very general theorem about sheaves $\cS$ of Abelian groups over a topological space $M$, for which a topology is introduced, in a certain way, on the cochain and cohomology groups $C^q(\fU,\cS)$, $H^q(\fU,\cS)$ of any open cover $\fU$ of $M$.
Let $\fV$ be another open cover, finer than $\fU$.
If the \v Cech differential has a certain openness property, then by Theorem 2.2 the refinement homomorphisms $H^q(\fU,\cS)\to H^q(\fV,\cS)$ are open.
To apply this result to cohesive sheaves over complex manifolds, all one has to do is to prove the required property of the \v Cech differential in (infinite dimensional) Stein manifolds.
This is the content of Theorem 4.6, whose proof rests on the vanishing theorem in \cite{LP}.

We hope to use the results of this paper to study group actions and relative duality in infinite dimensional manifolds, and analytic continuation in mapping spaces.

In the paper we will freely use basic sheaf theory and complex analysis, for which good references are \cite{Br,D2,M2,S3}.

\head 2.\ Cohomology groups of topologized sheaves \endhead

In this section $M$ is a Hausdorff space and
$\cS\to M$ is a sheaf of Abelian groups.
We will use standard notation of sheaf theory.
If $\fU=\{U_i\}_{i\in I}$ is an open cover of $M$ and $s=(i_0,\ldots,i_q)\in I^{q+1}$ is a $q$--simplex, $q\geq 0$, we write $U_s=U_{i_0}\cap\ldots\cap U_{i_q}$.
We also introduce a $(-1)$--simplex $s=\emptyset$, which constitutes $I^0$, and by $U_\emptyset$ we mean the whole space $M$.
Our goal is to endow the group
$$
C^q (\fU,\cS)=\prod_{s\in I^{q+1}}\ \Gamma(U_s,\cS),\quad q\geq -1,
$$
of (not necessarily alternating) cochains with a topology.
It would be natural to start by assuming the groups $\Gamma(U,\cS)$ are already endowed with a topology for open $U\subset M$, and give $C^q(\fU,\cS)$ the product topology.
However, with ulterior motives we take a more general track, and assume that $\cS$ is topologized in the following sense.

\definition{Definition 2.1}We say that the sheaf $\cS$ is topologized if for any Hausdorff space $U$ and any local homeomorphism $\pi\colon U\to M$ the group $\Gamma(U,\pi^{-1}\cS)$ is endowed with a topology, compatible with its group structure.
It is required that if $\rho\colon V\to U$ is another local homeomorphism, the induced map
$$
\rho^*\colon\Gamma (U,\pi^{-1}\cS)\to\Gamma(V,(\pi\rho)^{-1}\cS)
$$
be continuous.
\enddefinition

Given an open cover $\fU=\{U_i\}_{i\in I}$ of $M$ and $q\geq -1$, the space $\fU_q=\coprod_{s\in I^{q+1}} U_s$ admits a natural local homeomorphism $\pi_q$ into $M$, namely the one for which $\pi_q|U_s$ is the embedding $U_s\hookrightarrow M$.
We can then define a group topology on $C^q(\fU,\cS)$ as the image of the topology on $\Gamma(\fU_q,\pi^{-1}_q\cS)$ under the bijection
$$
\Gamma(\fU_q,\pi_q^{-1}\cS)\ni f\mapsto (f|U_s)_{s\in I^{q+1}}\in C^q(\fU,\cS).
$$
The \v Cech differential $\delta=\delta^q=\delta^q_{\fU}$
$$
\delta\colon C^q (\fU,\cS)\to C^{q+1} (\fU,\cS),\quad  q\geq 0,
$$
is a continuous homomorphism, and $Z^q(\fU,\cS)=\Ker \delta^q$ is a topological subgroup of $C^q(\fU,\cS)$.
This already defines the topology of $H^0(\fU,\cS)=Z^0 (\fU,\cS)$, and for $q\geq 1$ we endow the \v Cech group
$$
H^q(\fU,\cS)=Z^q (\fU,\cS)/\Im \delta^{q-1}
$$
with the quotient topology.
It will be advantageous to define $C^q(\fU,\cS)$ and $\delta^q$ for $q=-1$ as well: the former is simply $\Gamma(M,\cS)$, and
$$
\delta^{-1}\colon \Gamma(M,\cS)\ni f\mapsto (f|U_i)_{i\in I}\in C^0 (\fU,\cS).
$$ 

The cover $\fU$ induces a cover $\fU|W=\{U_i\cap W\}_{i\in I}$ of any $W\subset M$.
If $\fV=\{V_j\}_{j\in J}$ is another open cover of $M$, finer than $\fU$, any refinement map $\tau\colon J\to I$ defines continuous refinement homomorphisms
$$
C^q (\fU,\cS)\to C^q (\fV,\cS),\qquad H^q (\fU,\cS)\to H^q (\fV,\cS),
$$
that we denote $f\mapsto f|\fV$.
The map on cohomology is of course independent of the choice of $\tau$.
Write $\fV|\fU_p$ for the cover $\{\pi_p^{-1} V_j\}_{j\in J}$ of $\fU_p$.

\proclaim{Theorem 2.2}Suppose $\fU=\{U_i\}_{i\in I}$ and $\fV=\{V_j\}_{j\in J}$ are open covers of $M$, with $\fV$ finer than $\fU$, and for each $p,q\geq 0$ 
$$
\delta^{q-1}\colon C^{q-1} (\fV|\fU_p,\pi^{-1}_p\cS)\to Z^q (\fV|\fU_p,\pi_p^{-1}\cS)
$$
is an open map.
Then for $n\geq 0$ the refinement homomorphisms $H^n(\fU,\cS)\to H^n (\fV,\cS)$ are open.
More precisely, fixing any refinement map $\tau\colon J\to I$, the homomorphisms
$$
Z^n (\fU,\cS)\oplus C^{n-1} (\fV,\cS)\ni (f,g)\mapsto f|\fV+\delta_{\fV} g\in Z^n (\fV,\cS)
$$
are open for $n\geq 0$.
\endproclaim

For $q=0$ the assumption means that
$$
\delta^{-1}\colon\Gamma (\fU_p,\pi^{-1}_p\cS)\ni f\mapsto (f|\pi_p^{-1} V_j)_{j\in J}\in H^0 (\fV|\fU_p,\pi_p^{-1}\cS)\tag2.1
$$
is open, hence an isomorphism of topological groups.

Consider two complexes of Abelian topological groups $(A^\bullet,d_A)$ and $(B^\bullet,d_B)$ graded by $\bZ$, and a (continuous) homomorphism $\varphi\colon A^\bullet \to B^\bullet$.
Let $Z^\bullet (A^\bullet)=\Ker d_A$ and $Z^\bullet (B^\bullet)=\Ker d_B$.
We say that $\varphi$ is quasi-open if the homomorphism
$$
\varphi'\colon Z^\bullet (A^\bullet)\oplus B^\bullet\ni (a,b)\mapsto \varphi a-d_B b\in Z^\bullet (B^\bullet)\tag2.2
$$
is open, meaning that each $\varphi^{\prime n}\colon Z^n (A^\bullet)\oplus B^{n-1}\to Z^n (B^\bullet)$ is open.
Since $\varphi^{\prime n}$ is a homomorphism, this amounts to requiring that it map neighborhoods of $0\in Z^n (A^\bullet)\oplus B^{n-1}$ to neighborhoods of $0\in Z^n(B^\bullet)$.
In (2.2) one can recognize one component of the mapping cone of $\varphi$.
Clearly, a quasi-open $\varphi$ induces open maps $H^n(A^\bullet)\to H^n (B^\bullet)$ in cohomology.

The assumption of Theorem 2.2 is that $0^\bullet\to C^\bullet (\fV|\fU_p,\pi_p^{-1}\cS)$ is quasi-open, and its conclusion is that the refinement homomorphism $C^\bullet (\fU,\cS)\to C^\bullet (\fV,\cS)$ is quasi-open.

\proclaim{Proposition 2.3}Let $\varphi\colon A^\bullet\to B^\bullet$ and $\psi\colon B^\bullet\to C^\bullet$ be homomorphisms of complexes of Abelian topological groups.

\item{(a)}If $\varphi$ and $\psi$ are quasi-open, then so is $\psi\varphi$.
\item{(b)}If $\psi$ is a topological embedding and $0^\bullet\to C^\bullet/\psi B^\bullet$ is quasi-open, then $\psi$ is also quasi-open; and if in addition $\psi\varphi$ is quasi-open, then so is $\varphi$.
\endproclaim

\demo{Proof}We write $d$ for both differentials $d_B,d_C$. Now (a)\  follows from the formula
$$
(\psi\varphi)' (a,c)=\psi' (\varphi' (a,b), c-\psi b).
$$
As to (b), we can assume $B^\bullet\subset C^\bullet$ and $\psi$ is the inclusion map.
Let $\pi\colon C^\bullet\to C^\bullet/B^\bullet$ be the canonical projection, and $\zeta\colon 0^\bullet\to C^\bullet/B^\bullet$ the zero homomorphism.
If $z\in Z^n(C^\bullet)$ and $c\in C^{n-1}$ are such that $\zeta'(0,\pi c)=\pi z$, then $z+d c\in Z^n(B^\bullet)$ and $z=\psi' (z+d c,d c)$.
The map
$$
C^{n-1}\ni c\mapsto\zeta' (0,\pi c)\in Z^n (C^\bullet)
$$
being open, it follows that $\psi'$ is also open and $\psi$ is quasi-open.

Finally assume $\psi\varphi$ is also quasi-open.
Given zero neighborhoods $U_A\subset Z^n (A^\bullet)$ and $U_C\subset C^{n-1}$, we need to produce a zero neighborhood $U_B\subset Z^n (B^\bullet)$ contained in $\varphi'(U_A\times (B^{n-1}\cap U_C))$.
There are zero neighborhoods $U_1\subset C^{n-1}$ and $U_2\subset C^{n-2}$ satisfying $U_1+dU_2\subset U_C$; by shrinking $U_1$ we can arrange that
$$
\zeta' (\{0\}\times\pi U_2)\supset Z^{n-1} (C^\bullet/B^\bullet)\cap\pi U_1.
$$
Then $U_B=B^n\cap (\psi\varphi)' (U_A\times U_1)$ will do.
Indeed, it is open in $Z^n (B^\bullet)$, and any $z\in U_B$ can be written
$$
z=\varphi a-dc_1,\qquad a\in U_A,\ c_1\in U_1.
$$
Clearly $\pi dc_1=0$, so that $\pi c_1\in Z^{n-1} (C^\bullet/B^\bullet)\cap\pi U_1$, and $\pi c_1=\zeta' (0,\pi c_2)$ with some $c_2\in U_2$.
But then $b=c_1+dc_2\subset B^{n-1}\cap U_C$ and $z=\varphi a-db=\varphi'(a,b)$.
\enddemo

The rest of the section follows the spirit of \cite{S3, Chapter 1, \S4}.
Consider a double complex $K=\bigoplus_{p,q\geq 0} K^{pq}$ of Abelian topological groups with (continuous) differentials
$$
d'\colon K^{pq}\to K^{p+1,q},\qquad d''\colon K^{pq}\to K^{p,q+1},
$$
that satisfy $(d'+d'')^2=0$.
The grading $K^n=\bigoplus_{p+q=n} K^{pq}$ turns $K$ into an ordinary $\bZ$ graded complex (so $K^n=(0)$ if $n<0$), with differential $d=d'+d''$, whose cohomology groups $H^n(K^\bullet)$ will also be denoted $H^n (K)$.
Write $Z_I^{pq}(K)$, $Z_{II}^{pq}(K)$ for the kernel of $d'|K^{pq}$, resp.~$d''|K^{pq}$.
We also introduce subcomplexes $K_I, K_{II}\subset K$, where $K_I^{pq}=0$ if $q>0$ and $K_I^{p0}=Z_{II}^{p0}(K)$, while $K_{II}^{pq}=0$ if $p>0$ and $K_{II}^{0q}=Z_I^{0q}(K)$.
On $K_I$, resp.~$K_{II}$, the differential $d$ coincides with $d'$, resp.~$d''$, and $K_I^n=K_I^{n0}$, $K_{II}^n=K_{II}^{0n}$.

\proclaim{Proposition 2.4}Suppose that $d''\colon K^{pq}\to Z_{II}^{p,q+1}(K)$ is open for $p,q\geq 0$.
Then $0^\bullet\to K^\bullet/K_I^\bullet$, and hence the inclusion $K_I^\bullet\to K^\bullet$, are quasi-open. In particular, 
$H^n(K^\bullet_I)\to H^n(K^\bullet)$ is open for $n\ge 0$.
\endproclaim

\demo{Proof}By virtue of Proposition 2.3b it suffices to prove that $0^\bullet\to K^\bullet/K_I^\bullet$ is quasi-open; or that $0^\bullet\to K^\bullet$ is quasi-open, provided in addition to the original assumptions $K_I^\bullet=Z_{II}^{\bullet 0}(K)=0$ also holds.

Let $K_h=\bigoplus_{p\geq h} K^{pq}$.
These are subcomplexes of $K$, and $K_h/K_{h+1}$ is isomorphic to the double complex $L$ with $L^{pq}=0$ when $p\neq h$ and $L^{hq}=K^{hq}$; the differential on $L^\bullet$ agrees with $d''$.
From our assumptions $0^\bullet\to L^\bullet$ is quasi-open, hence by Proposition 2.3b so is the inclusion $K^\bullet_{h+1}\to K^\bullet_h$, and using Proposition 2.3a, so are all inclusions $K^\bullet_h\to K^\bullet_0=K^\bullet$.
Since $K_h^n=0$ when $h>n$, it now follows that $0^\bullet=\bigcap_h K_h^\bullet\to K^\bullet$ is also quasi-open, as claimed.
\enddemo

After these general considerations, let us return to the sheaf $\cS\to M$ of Theorem 2.2 and, at first, to two arbitrary open covers 
$\fU=\{U_i\}_{i\in I}$ and $\fV=\{V_j\}_{j\in J}$ of $M$.
Let $\frak W_{pq}=\coprod U_s\cap V_t$, the union over $s\in I^{p+1}$ and $t\in J^{q+1}$, and $\pi_{pq}\colon \frak W_{pq}\to M$ the local homeomorphism whose restriction to $U_s\cap V_t$ is the embedding $U_s\cap V_t\hookrightarrow M$.
Since $\cS$ will be fixed for the rest of the section, we will omit it from our notation, and write $\Gamma(U)$ for $\Gamma(U,\cS)$, $C^q(\fU)$ for $C^q(\fU,\cS)$, and so on.
Define a double complex $K=\bigoplus K^{pq}$,
$$
K^{pq}=\prod\Gamma (U_s\cap V_t),\qquad p,q\geq 0,
$$
the product over $s\in I^{p+1}$, $t\in J^{q+1}$.
Thus $K^{pq}$ is naturally identified with $\Gamma(\fW_{pq},\pi_{pq}^{-1}\cS)$, and we endow it with the topology induced by the topology on the latter.
There are continuous \v Cech--type differentials 
$$
d_{\fU}\colon K^{pq}\to K^{p+1,q}\quad\text{ and }\quad d_{\fV}\colon K^{pq}\to K^{p,q+1};
$$
for example if $f=(f_{st})\in K^{pq}$, the components of $d_{\fU}f=g$ are given by
$$
g_{i_0\ldots i_{p+1}t}=\sum^{p+1}_{k=0} (-1)^k f_{i_0\ldots\hat i_k\ldots i_{p+1}t}|U_{i_0\ldots i_{p+1}t}.
$$
Then $d'=d_{\fU}$ and $d''=((-1)^{q} d^{pq}_{\fV})_{pq}$ indeed turn $K$ into a double complex.
The map
$$
C^p(\fU)\ni (f_s)\mapsto (f_s|V_j)\in K_I^p=K_I^{p0}\tag2.3
$$
is a continuous isomorphism of groups.

\proclaim{Proposition 2.5}Suppose that for $p,q\geq 0$ 
$$
\delta^{q-1}\colon C^{q-1} (\fV|\fU_p)\to Z^p (\fV|\fU_p)\tag2.4
$$
are open.
Then (2.3) is an isomorphism of topological groups, and the inclusion $K^\bullet_I\to K^\bullet$ is quasi-open.
\endproclaim

\demo{Proof}The assumption for $q=0$ means (2.3) is open, hence a topological isomorphism.
As to the rest, the assumption for $q\geq 1$ means $d''\colon K^{pq}\to Z^{p,q+1}(K)$ are open for $p,q\geq 0$, and Proposition 2.4 completes the proof.
\enddemo

\proclaim{Proposition 2.6}If each component of $M$ is contained in some element of $\fU$,
then $\delta^{p-1}\colon C^{p-1}(\fU)\to Z^p (\fU)$ is open for $p\geq 0$.
\endproclaim

\demo{Proof}For each component $N$ of $M$ choose $k(N)\in I$ so that $N\subset U_{k(N)}$, and construct a continuous right inverse to $\delta^{p-1}$  by associating with $f=(f_t)\in Z^p(\fU)$ the cochain $g=(g_s)\in C^{p-1} (\fU)$, where $g_s|N=f_{k(N)s}|N$ for any component $N$ of $M$.
\enddemo

\demo{Proof of Theorem 2.2}Fixing a refinement map $\tau\colon J\to I$, let $\varphi\colon C^\bullet (\fU))\to C^\bullet (\fV)$ denote the refinement homomorphism.
Define continuous group isomorphisms $\alpha_I\colon C^\bullet (\fU)\to K^\bullet_I$ and $\alpha_{II}\colon C^\bullet (\fV)\to K^\bullet_{II}$ by
$$
\alpha_I (f_s)_{s\in I^{n+1}}=(f_s|V_j)_{s\in I^{n+1},j\in J},\ \ \alpha_{II} (g_t)_{t\in J^{n+1}}=(g_t|U_i)_{t\in J^{n+1},i\in I}.
$$
In fact, both are open, hence topological isomorphisms, $\alpha_I$ by Proposition 2.5, $\alpha_{II}$ by Proposition 2.6 (applied with $p=0$,
$M$ replaced by $\fV_n$ and $\fU$ by $\fU|\fV_n$).
Consider the compositions 
$$
\psi_I\colon C^\bullet (\fU)\overset{\alpha_I}\to\longrightarrow K^\bullet_I\hookrightarrow K^\bullet\text{ and }\psi_{II}\colon C^\bullet (\fV)\overset{\alpha_{II}}\to\longrightarrow K_{II}^\bullet\hookrightarrow K^\bullet.
$$
By Proposition 2.5 $\psi_I$ is quasi-open, and $\psi_{II}$ is a topological embedding.
Furthermore $\psi_I$ and $\psi_{II}\varphi$ induce the same map in cohomology.
More precisely, on top of p.~222 in \cite{S3}, Serre constructs a homomorphism $\gamma$ of degree $-1$,
$$
\gamma\colon Z^\bullet (\fU)\ni f\mapsto g^0-g^1+g^2-\ldots\in K^\bullet,
$$
which is clearly continuous in our set up, such that $\psi_{II}\varphi =\psi_I+d\gamma$.---Serre's notation is slightly different from ours.
He uses $\iota'$ both for our $\alpha_I$ and $\psi_I$, and similarly for $\iota''$; also, he denotes by $\tau$ the refinement homomorphism that we denoted $\varphi$.---For $f\in Z^n (\fU)$ and $k\in K^{n-1}$
$$
(\psi_{II}\varphi)' (f,k)=\psi_I (f)+d\gamma(f)-dk=\psi'_I (f,k-\gamma(f));
$$
since $\psi_I$ was quasi-open, so is $\psi_{II}\varphi$.
But by Proposition 2.4 (with the roles of I, II switched) and by Proposition 2.6 $0^\bullet\to K^\bullet/K^\bullet_{II}=K^\bullet/\psi_{II} C^\bullet (V)$ is also quasi-open, whence Proposition 2.3b implies $\varphi$ is quasi-open, as claimed.
\enddemo

\head 3.\ Cohesive sheaves\endhead

Here we quickly review the theory of cohesive sheaves, following \cite{LP}.
By a complex manifold we mean one that is modeled on a Banach space.
That is, a complex manifold $M$ is a Hausdorff space, sewn together from open subsets of Banach spaces, and the sewing maps are assumed to be holomorphic.
If $U\subset M$ is open and $E$ is a Banach space, the vector space of holomorphic functions $U\to E$ is denoted $\cO^E(U)$, and $\cO^E=\cO_M^E$ stands for the sheaf of $E$--valued holomorphic germs on $M$.
It will be convenient to distinguish between $\cO^E(M)$ and the space $\Gamma(M,\cO^E)$ of sections, although there is a bijection $e\mapsto\bold e$ between the two, where $\bold e(x)$ is defined as the germ of $e\in \cO^E(M)$ at $x\in M$.
The sheaves $\cO^E$ are called plain sheaves; they are modules over the sheaf $\cO=\cO^\bC$ of rings.
If $F$ is another Banach space, and $\Hom(E,F)$ stands the Banach space of continuous linear operators $E\to F$, any holomorphic map $M\to\Hom(E,F)$ induces an $\cO$-homomorphism $\cO^E\to\cO^F$ by composition.
Such homomorphisms are called plain, and form a module $\Hom_{\text{plain}}(\cO^E,\cO^F)$ over $\cO(M)$ (isomorphic to $\cO^{\Hom (E,F)}(M))$.
Germs of plain homomorphisms form a sheaf of $\cO$--modules, denoted {\bf Hom}$_{\text{plain}}(\cO^E,\cO^F)$, which in turn is isomorphic to $\cO^{\Hom (E,F)}$.

Suppose $X$ is a Banach space with a Schauder basis, and $\Omega\subset X$ is pseudoconvex, meaning $\Omega$ is open and $\Omega\cap Y$ is pseudoconvex in $Y$ for all finite dimensional subspaces $Y\subset X$.
A closed direct submanifold $N\subset\Omega$ will be called a Stein manifold, direct in the sense that $T_xN\subset T_x\Omega$ is a complemented subspace for every $x\in N$.
More generally, any complex manifold whose components are biholomorphic to such $N$ will be called Stein.
(This notion is related to what Patyi in \cite{P1} calls Banach--Stein manifold, but the two are not quite the same.)

\proclaim{Proposition 3.1}If open subsets $U_1,U_2$ of a complex manifold are Stein, then so is $U=U_1\cap U_2$.
\endproclaim

\demo{Proof}We can assume that there are Banach spaces $X_i$, pseudoconvex $\Omega_i\subset X_i$, closed direct submanifolds $N_i\subset\Omega_i$, and biholomorphisms $f_i\colon U_i\to N_i$, $i=1,2$.
It is then straightforward to check that $\Omega=\Omega_1\times\Omega_2\subset X_1\times X_2$ is pseudoconvex, the image of $f=(f_1|U,f_2|U)\colon U\to\Omega$ is a closed direct submanifold $N$ of $\Omega$, and $f\colon U\to N$ is a biholomorphism.
\enddemo

Fix a complex manifold $M$. Cohesive sheaves over $M$ are sheaves of $\cO$--modules with an extra structure and a special property. 
Below, sheaves of $\cO$--modules will be simply called $\cO$--modules. Given $\cO$--modules $\cA,\cB$, we write $\bHom_\cO(\cA,\cB)$ for the sheaf of $\cO$--homomorphisms between the two.

\definition{Definition 3.2}An analytic structure on an $\cO$--module $\cA$ is the choice, for each plain
sheaf $\cE$, of a submodule $\bHom(\cE,\cA)\subset\bHom_{\cO}(\cE,\cA)$, subject to
\item{(i)}if $\cE,\cF$ are plain sheaves and $\varphi\in\bHom_{\plain}(\cE,\cF)_z$ for some $z\in\Omega$,
then $\varphi^*\bHom(\cF,\cA)_z\subset\bHom(\cE,\cA)_z$; and
\item{(ii)}$\bHom(\cO,\cA)=\bHom_{\cO}(\cO,\cA)$.
\enddefinition

If $\cA$ is endowed with an analytic structure, we say that $\cA$ is an analytic sheaf.
The reader will realize that this is different from the traditional terminology, where ``analytic
sheaves'' and ``$\cO$--modules'' mean one and the same thing.

For example, one can endow a plain sheaf $\cG$ with an analytic structure by setting
$$
\bHom(\cE,\cG)=\bHom_{\plain}(\cE,\cG).
$$
We will always consider plain sheaves endowed with this analytic
structure.

An $\cO$--homomorphism $\varphi\colon\cA\to\cB$ of $\cO$--modules
induces a homomorphism
$$
\varphi_*\colon\bHom_{\cO} (\cE,\cA)\to\bHom_{\cO}(\cE,\cB)
$$
for $\cE$ plain.
When $\cA$, $\cB$ are analytic sheaves, we say that $\varphi$ is analytic if
$$
\varphi_*
\bHom(\cE,\cA)
\subset\bHom(\cE,\cB)
$$ 
for all plain sheaves $\cE$.
It is straightforward to check that if $\cA$ and $\cB$ themselves are plain sheaves, then $\varphi$ is
analytic precisely when it is plain.
We write $\Hom(\cA,\cB)$ for the $\cO(M)$--module of analytic homomorphisms $\cA\to\cB$ and $\bHom(\cA,\cB)$ for the sheaf of germs of analytic homomorphisms $\cA|U\to\cB|U$, with $U\subset M$ open.
Again, one easily checks that, when $\cA=\cE$ is plain, this new notation is consistent with the one already in use.
Further,
$$
\Hom(\cA,\cB)\approx\Gamma(M,\bHom(\cA,\cB)).
$$

\definition{Definition 3.3} Given an $\cO$-homomorphism
$\varphi:\cA\to\cB$ of $\cO$-modules,
any analytic structure on $\cB$ induces one on $\cA$ by the formula
$$
\bHom(\cE,\cA)=\varphi_*^{-1}\bHom(\cE,\cB).
$$
If $\varphi$ is an epimorphism, then any analytic structure on $\cA$ induces one on $\cB$ by the formula
$$
\bHom(\cE,\cB)=\varphi_*\bHom(\cE,\cA).
$$
\enddefinition
[LP, 3.4] explains this construction in the cases when $\varphi$ is the inclusion of a submodule
$\cA\subset\cB$ and when $\varphi$ is the projection on a quotient $\cB=\cA/\cC$.

\definition{Definition 3.4}A sequence $\cA\to\cB\to\cC$ of analytic sheaves and homomorphisms over
$M$ is said to be completely exact if for every plain sheaf $\cE$ and
every open
$U\subset M$, biholomorphic to a pseudoconvex set in a Banach space, the induced sequence 
$$
\Hom(\cE|U,\cA|U)\to\Hom(\cE|U,\cB|U)\to\Hom(\cE|U,\cC|U)\tag3.1
$$
is exact.
A general sequence of analytic homomorphisms is completely exact if every three--term subsequence is
completely exact.
An analytic homomorphism $\varphi\colon \cA\to \cB$ is a complete epimorphism if the sequence $\cA\overset\varphi\to\rightarrow\cB\to 0$ is completely exact.
\enddefinition

The above definition of complete exactness reduces to the definition in [LP, 4.1] when $M$ is an open subset of a Banach space.
Indeed, if [LP, 4.1] requires that (3.1) be exact just for $U$ pseudoconvex, in Definition 3.4 we are not requiring more, since an open subset of a Banach space biholomorphic to a pseudoconvex set is itself pseudodonvex; this follows from the characterization of pseudoconvexity in [M2, 3.75 Theorem (e) or (f)].

\definition{Definition 3.5}An infinite completely exact sequence
$$
\ldots\to\cF_2\to\cF_1\to\cS\to 0\
$$
of analytic homomorphisms is called a complete resolution of $\cS$ if each $\cF_j$ is plain.
\enddefinition

\definition{Definition 3.6}An analytic sheaf $\cS$ over $M$ is cohesive if each $z\in M$ has a neighborhood over which $\cS$ has a complete resolution.
\enddefinition

The simplest examples of cohesive sheaves are the plain sheaves, that have complete resolutions of form $\ldots\to 0\to 0\to\cE\to\cE\to 0$ and [LP] gives other examples of cohesive sheaves. By [L, Theorem 4.3] coherent sheaves over finite dimensional manifolds are also cohesive.
The main result of [LP] implies the following generalization of
Cartan's Theorems A and B:

\proclaim{Theorem 3.7}Let $\cS$ be a cohesive sheaf over a Stein manifold $M$.
Then
\itemitem{(a)}$\cS$ has a complete resolution $\ldots\to\cF_2\to\cF_1\to\cS\to 0$ over all of $M$;
\itemitem{(b)}$H^q (M,\cS)=0$ for $q\geq 1$; and more generally,
\itemitem{(c)}$H^q(M,\bHom(\cO_M^E,\cS))=0$ for $q\geq 1$ and Banach space $E$.
\endproclaim

\demo{Proof}It suffices to prove when $M$ is a direct submanifold of a pseudoconvex subset $\Omega$ of a Banach space with a Schauder basis.
By [P2], plurisubharmonic domination is possible in any pseudoconvex $\Omega'\subset\Omega$:\ given a locally bounded function $u\colon \Omega'\to\bR$, there is a continuous plurisubharmonic $v\colon\Omega'\to\bR$ such that $u\leq v$.
Therefore the assumptions of the Theorem in [LP, 12.2] are satisfied; the conclusion is that (a) and (b) indeed hold.
Part (c) follows the same way.
Let $\hat\cS$ denote the extension of $\cS$ to $\Omega$ by zero outside $M$.
By [LP, section 11, especially 11.7 and 11.10] and by [LP, 12.1], $\hat\cS\to\Omega$ is cohesive.
Since $\bHom(\cO_\Omega^E,\hat\cS)$ is the extension of $\bHom(\cO_M^E,\cS)$ by zero,
$$
H^q(M,\bHom (\cO_M^E,\cS))=H^q (\Omega,\bHom(\cO_\Omega^E,\hat\cS))=0,\quad q\geq 1,
$$
in view of [LP, 9.1].
\enddemo

We need to dwell on the completeness of the resolution $\cF_\bullet\to\cS\to 0$ of Theorem 3.7a.
In the setting when $M\subset\Omega$ is a direct submanifold as in the above proof, the proof in [LP], by reduction to [LP, 11.11], produces a sequence $\cF_\bullet\to\cS\to 0$ such that the induced sequence
$$
\Hom(\cO_U^E,\cF_\bullet|U)\to\Hom (\cO_U^E,\cS|U)\to 0\tag3.2
$$
is exact whenever $U=M\cap\hat U$ with $\hat U\subset\Omega$ pseudoconvex. Since the pair $(\Omega,M)$ is locally biholomorphic to a pair of a Banach space and its complemented subspace, it follows that each $x\in M$ has a neighborhood $V\subset M$ such that (3.2) is exact for any open $U\subset V$, biholomorphic to a pseudoconvex set; i.e., the resolution is completely exact over $V$.
But more is true:

\proclaim{Lemma 3.8}Suppose $M$ is a complex manifold, $\cS\to M$ is a cohesive sheaf, $\cF_n\to M$ are plain sheaves for $n\geq 1$, and
$$
\ldots\to\cF_2\overset{\varphi_2}\to\longrightarrow \cF_1\overset{\varphi_1}\to\longrightarrow\cS\overset{\varphi_0}\to\longrightarrow 0
$$
is a sequence of analytic homomorphisms, completely exact over some neighborhood of each $x\in M$.
Then for every Stein open $U\subset M$ and plain sheaf $\cE\to U$ the induced sequence
$$
\ldots\to\Hom (\cE,\cF_2|U)\to\Hom(\cE,\cF_1|U)\to\Hom(\cE,\cS|U)\to 0\tag3.3
$$
is exact.
\endproclaim

\demo{Proof}The analytic sheaves $\cK_n=\Ker\varphi_n|U$ fit in sequences
$$
0\to\cK_n\hookrightarrow\cF_n|U \overset{\varphi_n}\to\longrightarrow\cK_{n-1}\to 0,
$$
completely exact over some neighborhood of each $x\in M$.
By the Three Lemma [LP, 4.5] and by induction on $n$, the $\cK_n$ are cohesive; also
$$
0\to\bHom(\cE,\cK_n)\to\bHom (\cE,\cF_n|U)\to\bHom (\cE,\cK_{n-1})\to 0
$$
is exact.
Since in the associated exact sequence in cohomology
$$
0\to\Hom(\cE,\cK_n)\to\Hom (\cE,\cF_n|U)\to\Hom (\cE,\cK_{n-1})\to H^1 (U,\Hom (\cE,\cK_n))
$$
the last term is 0 by Theorem 3.7c, the sequence (3.3) is indeed exact.
\enddemo

Finally we have to discuss inverse images.
Let $\cS\to M$ be a cohesive sheaf, $U$ a Hausdorff space, and $\pi\colon U\to M$ a local homeomorphism.
We endow $U$ with the unique complex manifold structure that turns $\pi$ into a local biholomorphism.
Since locally $\cS$ and its inverse image $\pi^{-1}\cS$ are indistinguishable, the latter has canonical $\cO_U$--module and analytic structures; furthermore, with this analytic structure $\pi^{-1}\cS$ is cohesive.

\head 4.\ The topology on analytic cohomology groups\endhead

Suppose $M$ is a complex manifold modeled on complemented subspaces of Banach spaces with Schauder bases.
Such an $M$ will be called locally Stein, because the condition is equivalent to each $x\in M$ to have a Stein neighborhood.
Given a cohesive sheaf $\cS\to M$, we shall define a topology on the space $\Gamma(M,\cS)$ of sections, which then induces a topology on cochains $C^q(\fU,\cS)$ and on \v Cech groups $H^q(\fU,\cS)$ of open covers $\fU$ of $M$ and, when $M$ is paracompact, on $H^q(M,\cS)$.

First consider a plain sheaf $\cF=\cO^F$ over $M$ and a countable open cover $\fU=\{U_i\}_{i\in\bN}$ of $M$.
On the space
$$
\cO^F (\fU)=\{f\in\cO^F (M)\colon \|f\|_{U_i}=\sup_{U_i}\| f\|<\infty\text{ for }i\in\bN\}
$$
the seminorms $\|f\|_{U_i}$ define a Fr\'echet space structure.
Furthermore, every $f\in\cO^F(M)$ is contained in some $\cO^F(\fU)$, namely when $U_i=\{x\in M\colon \|f(x)\|< i\}$.
The locally convex direct limit of $\cO^F(\fU)$ is the finest locally convex topology on $\cO^F(M)$ for which all the inclusions $\cO^F(\fU)\hookrightarrow \cO^F(M)$ are continuous.
This topology is called the countable cover or $\tau_\delta$ topology, first introduced by Coeur\'e and Nachbin in [C,N].
Since $\tau_\delta$ is finer than the compact--open topology --- which is Hausdorff ---, it is also Hausdorff.
Neither $\cO^F(\fU)$ nor the Fr\'echet topology on it will change if we replace $U_i$ by $U'_i=\bigcup_{j\leq i} U_j$, so that in the above construction we can restrict to increasing countable covers $\fU$.
Since $\Gamma(M,\cO^F)$ is canonically identifiable with $\cO^F(M)$, we have also defined a topology on it, denoted again by $\tau_\delta$.
Clearly any plain homomorphism $\cO^E\to\cO^F$ induces a continuous linear map $\Gamma(M,\cO^E)\to\Gamma(M,\cO^F)$, and the same for maps $\pi^*\colon\Gamma (M,\cO^F)\to\Gamma (V,\cO^F)$ induced by a holomorphic map $\pi\colon V\to M$.

\definition{Definition 4.1}We call an analytic sheaf $\cS\to M$ separated if for  any open $U\subset M$, any plain sheaf $\cE\to U$, 
and any $\varphi\in\Hom (\cE,\cS|U)$
\itemitem{(i)}$\Ker\varphi_*\subset\Gamma (U,\cE)$ is closed;
\itemitem{(ii)}if $\varphi\be=0$ for all constant germs $\be\in \cE$, then $\varphi=0$.
\enddefinition

This is clearly a local property of $\cS$.
The notion will be important for cohesive sheaves.
Not all cohesive sheaves are separated, as we shall see in section 6 (although we do not know if (ii) can fail for a cohesive $\cS$).
Still, many are:

\proclaim{Lemma 4.2}(a)\ Plain sheaves are separated.
(b)\ Analytic subsheaves of separated sheaves are separated.
(c)\ If $N\subset M$ is a direct submanifold, $F$ is a Banach space, $k\in\bN$, and $\cJ\subset\cO_M^F=\cO^F$ is the subsheaf of germs vanishing on $N$ to order $k$, then $\cO^F\!/\!\cJ$ is separated.
\endproclaim

\demo{Proof}(a)\ Consider a plain sheaf $\cO^F\to M$.
Any $\varphi\in\Hom(\cO^E|U,\cO^F|U)$ is induced by some $\Phi\in\cO^{\Hom(E,F)}(U)$, and first we need to show that $\{e\in \cO^E(U)\colon\Phi e=0\}$ is a closed subspace of $\cO^E(U)$; which is clear since it is closed even in the topology of pointwise convergence.
Second we need to show that if $\Phi e=0$ for all constant functions $e\colon U\to E$, then $\Phi\equiv 0$, which again is obvious.

(b)\ If $\cA$ is an analytic sheaf and $\cB\subset\cA$ an analytic subsheaf, the claim follows from
$$
\Hom(\cE|U,\cB|U)\subset\Hom(\cE|U,\cA|U).
$$

(c)\ It suffices to verify the requirements of Definition 4.1 for $U=M$, and we might as well assume $M$ is Stein.
By [LP, 10.3] $\cJ$ is cohesive, and by the definition of the induced analytic structure on $\cJ$, for plain $\cO^E\to M$
$$
0\to\bHom (\cO^E,\cJ)\to\bHom (\cO^E,\cO^F)\to\bHom (\cO^E,\cO^F\!/\!\cJ)\to 0
$$
is exact.
Hence in the associated long exact sequence
$$
\ldots\to\Hom (\cO^E,\cO^F)\overset p\to\longrightarrow \Hom (\cO^E,\cO^F\!/\!\cJ)\to H^1 (M,\bHom (\cO^E,\cJ))\to\ldots
$$
the last term is 0 by Theorem 3.7, and so $p$ is onto.
Given $\varphi\in\Hom (\cO^E,\cO^F\!/\!\cJ)$, choose $\psi\in\Hom (\cO^E,\cO^F)$ so that $p\psi=\varphi$.
This $\psi$ is induced by a holomorphic $\Psi\colon M\to\Hom (E,F)$, and in the identification $\Gamma(M,\cO^E)\leftrightarrow \cO^E(M)$, the kernel of $\varphi_*$ corresponds to
$$
J=\{e\in\cO^E (M)\colon\Psi e\text{ vanishes on }N\text{ to order }k\}.
$$
Clearly $J\subset \cO^E(M)$ is closed in the compact--open topology, hence in the finer $\tau_\delta$ topology as well.
Further, if $J$ contains all constant functions $e\colon M\to E$, then $\Psi$ itself vanishes to order $k$ on $N$, i.e., $\psi\cO^E\subset\cJ$ and $\varphi=0$ as claimed.
\enddemo

Now consider an arbitrary $\cO$--module $\cS\to M$.
Any $\cO$--homomorphism $\varphi\colon\cO^F\to\cS$
induces a topology $\varphi_*\tau_\delta$ on $\Gamma(M,\cS)$.
This is the finest locally convex topology for which $\varphi_*\colon\Gamma (M,\cO^F)\to\Gamma (M,\cS)$ is continuous, when $\Gamma(M,\cO^F)$ is endowed with the $\tau_\delta$ topology.
If $\psi\in\Hom_{\plain} (\cO^E,\cO^F)$ and $\var=\varphi\psi\colon\cO^E\to\cS$, then the topology $\varphi_*\tau_\delta$ is coarser than $\var_*\tau_\delta$.
We apply this when $M$ is Stein, $\cS$ is cohesive, and $\varphi\colon\cO^F\to\cS$ is a complete epimorphism (which exists by Theorem 3.7a).
If $\var\in\Hom (\cO^E,\cS)$, by Definition 3.4 there is a $\psi\in\Hom_{\plain}(\cO^E,\cO^F)$ such that $\var=\varphi \psi$.
Therefore $\varphi_*\tau_\delta$ is coarser than $\var_*\tau_\delta$, indeed it is the coarsest topology induced from $\tau_\delta$ by an analytic homomorphism $\cO^E\to\cS$.
In particular, $\varphi_*\tau_\delta$ is independent of the choice of the complete epimorphism $\varphi\colon\cO^F\to\cS$, and for brevity we call it again the $\tau_\delta$ topology on $\Gamma(M,\cS)$.
Remember that for the time being, $M$ is Stein.
When $\cS$ itself is plain, the new definition agrees with the old one, since we can take as $\varphi$ the complete epimorphism $\id\colon\cS\to\cS$.

The following is straightforward:

\def\cT{\Cal T}
\proclaim{Proposition 4.3}Let $M$ and $V$ be Stein manifolds, $\cS,\cT$ cohesive sheaves over $M$, and $\pi\colon V\to M$ a local biholomorphism (so that $\pi^{-1}\cS\to V$ is cohesive, see the last paragraph in section 3).
Then any analytic homomorphism $\cS\to\cT$ induces a continuous map $\Gamma(M,\cS)\to \Gamma(M,\cT)$, and 
pull back $\pi^*\colon\Gamma(M,\cS)\to\Gamma (V,\pi^{-1}\cS)$ is continuous.
\endproclaim

Although the definitions already made suffice to formulate our main result, for the sake of completeness we now define the topology on $\Gamma(M,\cS)$ for locally Stein $M$. 

\definition{Definition 4.4}If $\cS$ is a cohesive sheaf over a locally Stein manifold $M$, the $\tau^\delta$ topology on $\Gamma(M,\cS)$ is the coarsest locally convex topology for which restriction $\Gamma(M,\cS)\to\Gamma(U,\cS)$ is continuous when $\Gamma(U,\cS)$ is endowed with the $\tau_\delta$ topology, for every Stein open $U\subset M$.
\enddefinition

When $M$ itself is Stein, it follows from Proposition 4.3 that the $\tau_\delta$ and $\tau^\delta$ topologies on $\Gamma(M,\cS)$ agree.
For general locally Stein $M$, on $\Gamma(M,\cO^F)$ $\tau_\delta$ is finer than $\tau^\delta$, but the two may differ as far as we can tell.

Comparing Definition 4.4 and Proposition 4.3 we see that the latter generalizes to local biholomorphisms $\pi\colon V\to M$ of locally Stein manifolds:\ for a cohesive $\cS\to M$ the induced map $\pi^*\colon\Gamma (M,\cS)\to\Gamma (V,\pi^{-1}\cS)$ is continuous in 
$\tau^\delta$. The $\tau^\delta$ topology makes $\cS\to M$ a topologized sheaf in the sense of Definition 2.1, and defines a locally convex topology on the cochain and \v Cech groups $C^q(\fU,\cS)$, $H^q (\fU,\cS)$ for an open cover $\fU$ of $M$, as explained in section 2.
It also defines a topology on the direct limit $\check H^q (M,\cS)=\varinjlim_{\fU} H^q (\fU,\cS)$, the finest locally convex topology on 
$\check H^q (M,\cS)$ for which the canonical maps
$$
H^q (\fU,\cS)\to\check H^q(M,\cS)\tag4.1
$$
are continuous.
We refer to these topologies again as the $\tau^\delta$ topologies.
Of course, when $M$ is paracompact, the sheaf cohomology group $H^q(M,\cS)$ and $\check H^q(M,\cS)$ are canonically isomorphic, so one can speak of the $\tau^\delta$ topology on the former as well.
Our main result is

\proclaim{Theorem 4.5}Suppose $M$ is a locally Stein manifold, $\cS\to M$ a separated cohesive sheaf, $\fU$ a cover of $M$ by Stein open sets, and $q=0,1,2,\ldots$.
Then the canonical map (4.1) is an isomorphism of topological vector spaces.
\endproclaim

In section 6 we will show that the claim does not hold for all cohesive sheaves.---
Theorem 4.5 can be quickly derived from the following result, which we will prove in section 5:

\proclaim{Theorem 4.6}Let $M$ be a Stein manifold, $\cS\to M$ a separated cohesive sheaf, $\fV=\{V_j\}_{j\in J}$ a cover of $M$ by Stein open subsets, and $q=0,1,2,\ldots$
Then
$$
\delta\colon C^{q-1} (\fV,\cS)\to Z^q (\fV,\cS)
$$
is open (in the $\tau^\delta$ topology).
\endproclaim

Assuming Theorem 4.6, the following is an immediate consequence, in view of Theorem 2.2 and Proposition 3.1:

\proclaim{Proposition 4.7}Let, $M,\cS,\fU=\{U_i\}_{i\in I}$, $q$ be as in Theorem 4.5, and let $\fV=\{V_j\}_{j\in J}$ be another cover of $M$ by Stein open subsets, finer than $\fU$.
Fixing a refinement map $J\to I$, the induced map
$$
Z^q (\fU,\cS)\oplus C^{q-1} (\fV,\cS)\ni (f,g)\mapsto f|\fV+\delta g\in Z^q (\fV,\cS)
$$
is open, and therefore so is the refinement homomorphism $H^q (\fU,\cS)\to H^q (\fV,\cS)$.
\endproclaim

\demo{Proof of Theorem 4.5}Let $\fV=\{V_j\}_{j\in J}$ be another cover of $M$ by Stein open subsets, finer than $\fU$.
Since $U_s$ and $V_t$ for $p\geq 1$ and $s\in I^p$, $t\in J^p$ are Stein (Proposition 3.1), by Theorem 3.7 and by Leray's theorem,
$$
H^q (\fV|U_s,\cS)\approx H^q (U_s,\cS)=0,\quad\text{for }q\geq 1.
$$
Therefore, e.g.~by [S3, Proposition 5], the refinement homomorphism $H^q (\fU,\cS)\to H^q (\fV,\cS)$ is bijective and clearly continuous.
It is also open by Proposition 4.7, hence an isomorphism of topological vector spaces.
Since $\check H^q (M,\cS)$ can be obtained as the direct limit of $H^q(\fV,\cS)$ with covers $\fV$ by Stein open subsets, the canonical map $H^q(\fU,\cS)\to\check H^q (M,\cS)$ is indeed an isomorphism of topological vector spaces.
\enddemo

\head 5.\ The proof of Theorem 4.6\endhead

It will be useful to represent the countable cover topology on $\Gamma(M,\cS)$ somewhat differently than in the original construction.
Suppose $M$ is an arbitrary complex manifold.
Let us call a locally bounded function $w\colon M\to (0,\infty)$ a weight, and write $W=W(M)$ for the set of weights.
This is a directed set for the partial order $u<v$ meaning $u(x)< v(x)$ for all $x\in M$.
If $F$ is a Banach space, set
$$
\cO^F(w)=\{f\in\cO^F (M)\colon \|f\|_w=\sup_{x\in M} \|f(x)\|/w(x) < \infty\},
$$
so that $(\cO^F(w),\ \|\ \|_w)$ is a Banach space.
Since $f\in\cO^F (\|f\|)$ for $f\in\cO^F(M)$, and $\cO^F(u)\subset\cO^F(v)$ when $u< v$, as a vector space $\cO^F(M)$ is the direct limit of its subspaces $\cO^F(w)$.

\proclaim{Proposition 5.1}$\underset{w}\to\varinjlim (\cO^F(w), \|\ \|_w)=(\cO^F(M),\tau_\delta)$ as locally convex spaces.
\endproclaim

\demo{Proof}Given $w\in W$, consider the cover $\fU$ of $M$ by $U_i=\{x\in M\colon w(x) < i\}$, $i\in\bN$.
As the inclusion $\cO^F(w)\hookrightarrow \cO^F (\fU)$ is continuous, so are $\cO^F(w)\hookrightarrow \cO^F (M)$ and the identity map $\varinjlim \cO^F(w)\to\cO^F (M)$.
Conversely, start with an open cover $\fU=\{U_i\}_{i\in\bN}$.
We claim that $\cO^F(\fU)\hookrightarrow \varinjlim \cO^F(w)$ is continuous, i.e., any neighborhood $\cU$ of 0 in the latter space contains a neighborhood of $0\in\cO^F (\fU)$.
Now $\cU$ contains a set of form
$$
\cU_\var=\bigcup_{w\in W} \{f\in\cO^F (M)\colon \|f\|_w <\var (w)\},
$$
where $\var\colon W\to (0,\infty)$.
Suppose there is no neighborhood of $0\in\cO^F(\fU)$ that is contained in $\cU$; then for each $n\in\bN$ there is an $f_n\in\cO^F (M)\backslash\cU$ such that
$$
\sup_{U_i} \|f_n\| < 1/n,\quad\text{for }i=1,\ldots,n.
$$
If $\overline w(x)=1+\sup_n n\|f_n(x)\|$, then $\|f_n\|_{\overline w}< 1/n$, and therefore $f_n\in\cU_\var\subset\cU$ when $n>1/\var (\overline w)$, after all.
This contradiction proves that $\cO^F(\fU)\hookrightarrow\varinjlim\cO^F(w)$ is continuous, hence so is the identity map $\cO^F(M)\to\varinjlim \cO^F(w)$; the proof is complete.
\enddemo

Under the identification $\cO^F(M)\leftrightarrow\Gamma (M,\cO^F)$ to the spaces $\cO^F(w)$ correspond Banach subspaces $\Gamma(w,\cO^F)\subset \Gamma(M,\cO^F)$, for which the norm will also be denoted $\|\ \|_w$.
Suppose $\cS\to M$ is a cohesive sheaf over a Stein manifold, and fix a complete epimorphism $\varphi\colon\cO^F\to\cS$.
If $w\in W(M)$, we denote by $\Gamma(w,\cS)$ the image of $\Gamma(w,\cO^F)$ under $\varphi_*\colon\Gamma(M,\cO^F)\to\Gamma(M,\cS)$, and by $\|\ \|_w$ the seminorm on $\Gamma(w,\cS)$
$$
\|\sigma\|_w=\inf \{\|f\|_w\colon\varphi f=\sigma\},\qquad \sigma\in\Gamma(w,\cS).
$$
This seminorm induces the quotient topology on $\Gamma(w,\cS)$, Hausdorff when $\cS$ is separated.
Similarly, let $\fU=\{U_i\}_{i\in I}$ be a Stein cover of $M$ and $q\geq -1$.
As in section 2, we consider $\fU_q=\coprod_{s\in I^{q+1}} U_s$ and its natural map $\pi_q\colon\fU_q\to M$.
In our case $\fU_q$ is a Stein manifold and $\pi_q$ a local biholomorphism; $\pi_q^{-1}\cS$ is cohesive and $\varphi$ induces a complete epimorphism $\varphi_q\colon\cO^F_{\fU_q}\to\pi_q^{-1}\cS$.
The construction above defines, given a weight $w\in W(\fU_q)$, the space $\Gamma(w,\pi_q^{-1}\cS)$ and the seminorm $\|\ \|_w$ on it.
Identifying $\Gamma(\fU_q,\pi^{-1}_q\cS)$ with $C^q(\fU,\cS)$, we obtain the subspace $C^q(w,\cS)\subset C^q (\fU,\cS)$ corresponding to $\Gamma(w,\pi_q^{-1}\cS)$, and a seminorm, again denoted $\|\ \|_w$, on it.
Proposition 5.1 implies $C^q(\fU,\cS)=\varinjlim C^q (w,\cS)$ as locally convex spaces.

\demo{Proof of Theorem 4.6}Fix a complete epimorphism $\varphi\colon\cO^F\to\cS$.
Since $\cS$ is separated, $(C^q (w,\cS),\ \|\ \|_w)$ is a Banach space when $q\geq -1$ and $w\in W(\fV_q)$; further
$$
Z^q (w,\cS)=C^q (w,\cS)\cap Z^q (\fV,\cS)
$$
is a closed subspace and
$$
Z^q (\fV,\cS)=\underset w\to\varinjlim  (Z^q (w,\cS), \|\ \|_w)
$$
as locally convex spaces.
Now in the direct limit $C^q(\fV,\cS)$ a neighborhood basis of 0 is formed by the convex hull of sets of form
$$
\cU_\var=\bigcup_{w\in W(\fV_q)} \{\sigma\in C^q (w,\cS)\colon \|\sigma\|_w <\var(w)\},
$$
where $\var\colon W(\fV_q)\to (0,\infty)$; and similarly for $Z^q(\fV,\cS)$.
Hence to prove the theorem it will suffice to produce for every $w\in W(\fV_q)$ a $w'\in W(\fV_{q-1})$ with the property that given $\sigma\in Z^q(w,\cS)$, there is a $\sigma'\in C^{q-1} (w',\cS)$ such that
$$
\delta\sigma'=\sigma\quad\text{and}\quad \|\sigma'\|_{w'}\leq \|\sigma\|_w.\tag5.1
$$

Given $w\in W(\fV_q)$, consider the closed subspace
$$
E=\{e\in C^q (w,\cO^F)\colon\delta\varphi e=0\}\subset C^q (w,\cO^F),
$$
with the inherited norm $\|e\|_E=\|e\|_w$.
There is a tautological cochain $\tau\in C^q(\fV,\bHom (\cO^E,\cO^F))$ obtained as follows.
If $s\in J^{q+1}$, $x\in V_s$, and $e=(e_t)_{t\in J^{q+1}}\in E$, set
$$
\tau_s(x) (e)=e_s(x).
$$
This defines a map $\tau_s\colon V_s\to\Hom(E,F)$ that is holomorphic, as one easily checks.
The tautological cochain $\tau$ then has components $\tau_s$ (or rather, the plain homomorphisms induced by $\tau_s$).
That is, if we denote by $\be \in\Gamma (M,\cO^E)$ the section corresponding to the constant function $\equiv e$, then
$$
\tau\be =e.\tag5.2
$$
Let $\cI\in C^q (\fV,\bHom(\cO^E,\cS))$ denote the image of $\tau$ under $\varphi$.

The point is that $\delta\cI=0$.
To see this, let $e\in E$.
Using (5.2)
$$
(\delta\cI)\be=\delta(\cI \be)=\delta(\varphi\tau \be)=\delta(\varphi e)=0.
$$
Since the components of $\delta\cI$ are sections of $\bHom(\cO^E,\cS)$, and $\cS$ is separated, it follows that indeed $\delta\cI=0$.

As $H^q(\fV,\bHom(\cO^E,\cS))=0$ for $q\geq 1$ follows from Theorem 3.7, there is $\chi\in C^{q-1} (\fV,\bHom (\cO^E,\cS))$ with $\delta\chi=\cI$; and this is also true for $q=0$ by the very definition of a sheaf.
Since $\varphi$ is a complete epimorphism, there is a $\kappa\in C^{q-1} (\fV,\bHom (\cO^E,\cO^F))$ such that $\varphi\kappa=\chi$.
The components of $\kappa$ define a holomorphic function $h\colon \fV_{q-1}\to\Hom (E,F)$, and we claim that $w'=2\| h\|\in W(\fV_{q-1})$ will do.
(Here $\|h(x)\|$ is the operator norm of $h(x)\in\Hom (E,F)$.)
Indeed, suppose $\sigma\in Z^q (w,\cS)$.
There is an $e\in C^q (w,\cO^F)$ such that $\varphi e=\sigma$ and $\|e\|_w\leq 2\|\sigma\|_w$; clearly $e\in E$.
As before, denoting by $\be\in \Gamma(M,\cO^E)$ the section corresponding to the constant function $\equiv e$, letting $\sigma'=\chi\be \in C^{q-1} (\fV,\cS)$, and using (5.2)
$$
\delta\sigma'=\delta\chi \be=\cI\be=\varphi\tau\be=\varphi e=\sigma.
$$
Further, $\sigma'=\varphi\kappa\be$.
Now $\kappa\be \in C^{q-1} (\fV,\cO^F)$ corresponds to the holomorphic function $h e\colon \fV_{q-1}\to F$.
Therefore
$$
\|\sigma'\|_{w'}\leq \|\kappa\be \|_{w'}=\|he\|_{w'}
=\sup_{x\in\fV_{q-1}}\frac{\|h(x) e\|_F}{w' (x)}\leq\sup_{x\in\fV_{q-1}} \frac{\|h(x)\| \|e\|_E}{w'(x)}\leq \|\sigma\|_{w}
$$
and the proof is complete.
\enddemo

\head 6.\ An example\endhead

In this section we give an example of a cohesive (but nonseparated) sheaf $\cS\to M$ and a Stein cover $\fV$ of $M$ such that $H^q(\fV,\cS)\to \check H^q(M,\cS)=H^q(M,\cS)$ is not a homeomorphism.
In fact, $M$ will be $\bC$, $\fV$ will consist of the single
set $\bC$, and $q=0$.
We start with two general results.

\proclaim{Lemma 6.1}Let $X,E,F$ be Banach spaces, $B\subset X$ a ball, $\xi\in\partial B$, and $P\colon B\to\Hom (E,F)$ a holomorphic function.
If for every $e\in E$ the function $Pe\colon B\to F$ analytically continues across $\xi$, then $P$ itself continues across $\xi$.
\endproclaim

A similar result for one dimensional $X$ was already proved in [LP, page 464].

\demo{Proof}We can assume $B$ is the unit ball.
Let us write $\partial^k$ for the $k$'th iterated derivative in the direction $\xi$.
Given $\var, p>0$, the set
$$
E_{p\var}=\{e\in E\colon\text{ if } \|\eta-\xi/2\|_X <\var\text{ and }k\in\bN,\ \text{then}\ \|\partial^k P(\eta) e\|_F\leq p k! (1/2+\var)^{-k}\}
$$
is closed.
Our assumptions imply $\bigcup_{p,\var>0} E_{p\var}=E$, whence one of the $E_{p\var}$ has an interior point, and then $0\in\text{int } E_{2p,\var}$.
It follows that for some $q>0$
$$
\|\partial^k P(\eta)\|_{\Hom(E,F)}\leq q k! (1/2+\var)^{-k},\text{ whenever } \|\eta-\xi/2\|_X <\var.\tag6.1
$$
Now let $l\colon X\to\bC$ be a linear form, $l(\xi)=1/2$.
The estimate (6.1) implies that the series
$$
Q(x)=\sum^\infty_{k=0}\partial^k P(x-l(x)\xi) l(x)^k\!/\!k!
$$
defines a function $Q$ holomorphic in a neighborhood of $\xi$, which provides the required analytic continuation of $P$.
\enddemo

\proclaim{Lemma 6.2}Let $M$ be a locally Stein manifold, $E,F$ Banach spaces, and $\Phi\in\cO^{\Hom(E,F)}(M)$. If $\Phi(x)$ is injective for $x$ in a dense open $D\subset M$,
then the homomorphism $\varphi\colon\cO^E\to\cO^F$ induced by $\Phi$ is an analytic isomorphism between $\cO^E$ and $\varphi\cO^E$, the latter endowed with the analytic structure inherited as a subsheaf of $\cO^F$.
\endproclaim

\demo{Proof}Clearly $\varphi$ is an $\cO$--isomorphism $\cO^E\to\varphi \cO^E$.
To prove the lemma, it will suffice to show that if $G$ is a Banach space, $U\subset M$ is open, and $\Psi\in\cO^{\Hom(G,F)}(U)$ is such that for every $\gamma\in\cO^G(U)$ the germs of $\Psi\gamma$ are in $\varphi\cO^E$, then there is a $P\in\cO^{\Hom(G,E)}(U)$ such that $\Psi=\Phi P$.

To verify this, for $x\in D\cap U$ define a linear operator $\Pi (x)\colon G\to E$ by letting
$$
\Pi(x) g=e,\quad\text{ if }\quad \Psi (x) g=\Phi (x) e.\tag6.2
$$
One checks that $\Pi(x)$ is closed, hence continuous by the closed graph theorem.
Further, given $g\in G$, the germs of $\Psi g$ are in $\varphi\cO^E$, so for any $\xi\in U$ there are a neighborhood $V\subset U$ of $\xi$ and $\var\in\cO^E (V)$ such that $\Psi g=\Phi\var$ on $V$.
If $\xi\in D\cap U$, (6.2) implies $\Pi g=\var$ is holomorphic near $\xi$; hence $\Pi g$ is holomorphic on $D\cap U$.
This in turn implies $\Pi\colon D\cap U\to\Hom(G,E)$ is holomorphic, see [M, Exercise 8E], whose solution rests on Cauchy's formula and the principle of uniform boundedness.

We next show that $\Pi$ extends to a $P\in\cO^{\Hom(G,E)}(U)$.
For this we can assume $U$ is an open subset of a Banach space.
Consider pairs $(D',\Pi')$ consisting of an open $D'\subset U$ containing $D$ and $\Pi'\in\cO^{\Hom(G,E)}(D')$ extending $\Pi$.
Let us write $(D',\Pi')< (D'',\Pi'')$ if $D'\subset D''$ and $\Pi'=\Pi''|D'$.
Zorn's lemma, applied to this partial order, gives a maximal element $(D_0,P)$.
As $D$ is dense, (6.2) implies $\Psi=\Phi P$ on $D_0$.
But $D_0$ must be $U$ itself.
Otherwise there would be a ball $B\subset D_0$ with a point $\xi\in U\backslash D_0$ on its boundary.
Take any $g\in G$.
As above, there are a neighborhood $V\subset U$ of $\xi$ and $\var\in\cO^E(V)$ such that $\Psi g=\Phi\var$ on $V$.
By (6.2) this implies $Pg=\var$ on $D\cap V$, and so $Pg$ analytically continues across $\xi$.
In view of Lemma 6.1 therefore $P$ itself continues to a neighborhood $W$ of $\xi$.
Denote this continuation $P'\in\cO^{\Hom(G,E)}(W)$.
As $\Psi=\Phi P$ on $D_0$, also $\Psi=\Phi P'$ on $W$.
But $\Phi$ is injective on $D\cap U$, therefore $P'=P$ on $D\cap W$, and by density, on $D_0\cap W$.
Hence $P$ analytically continues to $D_0\cup W$, contradicting the maximality of $(D_0,P)$.

Therefore indeed $P\in\cO^{\Hom(G,E)}(U)$, and the proof is complete.
\enddemo

\proclaim{Theorem 6.3}There is a cohesive sheaf $\cS\to\bC$ such that the canonical map $\Gamma(\bC,\cS)\to H^0(\bC,\cS)$ is not a homeomorphism.
\endproclaim

This is equivalent to saying that the canonical map $H^0 (\fU,\cS)\to H^0 (\bC,\cS)$ is not a homeomorphism when $\fU=\{\bC\}$.

\demo{Proof}
If $M$ is a finite dimensional manifold, $\fU$ its cover by relatively compact open subsets, and $\fV$ is a finer cover, then not only the (bijective) refinement homomorphism $\cO^E(\fU)\to\cO^E(\fV)$ is continuous, but its inverse, too.
It follows that the canonical map $\cO^E(\fU)\to\cO^E(M)$ is a homeomorphism, and so $\cO^E(M)$ and $\Gamma(M,\cO^E)$ are Fr\'echet spaces.
Hence given a complete epimorphism $\pi\colon\cO^E\to\cS$ on a cohesive sheaf, a sequence $\sigma_n\in\Gamma (M,\cS)$ converges if and only if it is the image under $\pi$ of a convergent sequence in $\Gamma(M,\cO^E)$.

Now let $E=l^p$ with some $p\in [1,\infty]$, and consider $\Phi\in \cO^{\Hom(E,E)}(\bC)$ given by
$$
\Phi(\zeta)(x)=(x_n\zeta^{n+1}/n^{n+1})^\infty_{n=1},\qquad\zeta\in\bC,\ x=(x_n)\in E.
$$
Let $\varphi\colon\cO^E\to\cO^E$ be the plain homomorphism induced by $\Phi$, $\cS=\cO^E/\varphi\cO^E$, and $\pi\colon\cO^E\to\cS$ the canonical projection.
Lemma 6.2 implies that for any plain sheaf $\cO^G$ the sequence
$$
0\to\bHom(\cO^G,\cO^E)\overset{\varphi_*}\to\rightarrow\bHom (\cO^G,\cO^E)\overset{\pi_*}\to\rightarrow
\bHom(\cO^G,\cS)\to 0
$$
is exact, and it follows that
$$
0\to\cO^E\overset\varphi\to\rightarrow \cO^E\overset\pi\to\rightarrow\cS\to 0
$$
is completely exact ([L, Theorem 2.7]). Hence $\cS$ is cohesive.

The canonical map $\Gamma(\bC,\cS)\to H^0 (\bC,\cS)$ is of course continuous, so the point of the theorem is that its inverse is discontinuous.
We cover $\bC$ by two sets
$$
U\subset \{\zeta\in\bC\colon |\zeta| < 1\}\quad\text{ and }\quad V\subset\bC\backslash \{0\},
$$
and show that the inverse of the map
$$
\Gamma(\bC,\cS)\ni s\mapsto (s|U, s|V)\in H^0 (\{U,V\},\cS)\tag6.3
$$
is discontinuous.
This will then prove the theorem, since the canonical map\newline
$H^0(\{U,V\},\cS)\to H^0 (\bC,\cS)$ is continuous.

For $n=1,2,\ldots$ consider the holomorphic function $e_n\colon\bC\to E$
$$
e_n(\zeta)=(0,\ldots,\zeta^n/n,0,\ldots),
$$
the only nonzero entry appearing at the $n$'th spot.
Clearly $e_n|U\to 0$ uniformly.
Hence with the corresponding sections $\be_n\in\Gamma(\bC,\cO^E)$ then $\be_n|U\to 0$ in $\Gamma(U,\cO^E)$ and $\pi\be_n|U\to 0$ in $\Gamma(U,\cS)$.
Also $\pi\be_n|V=0$, because
$$
e_n=\Phi f_n,\text{ where }f_n(\zeta)=(0,\ldots,n^n/\zeta,0,\ldots).
$$
However, $\pi\be_n\not\rightarrow 0$ in $\Gamma(\bC,\cS)$ (and therefore the inverse of (6.3) is indeed discontinuous).

For suppose $g_n\in\cO^E(\bC)$ project to $\pi\be_n$, that is, $\pi\bg_n=\pi\be_n$.
This means $\bg_n-\be_n\in\varphi\cO^E$, or
$$
g_n-e_n=\Phi h_n\quad\text{with some}\quad h_n\in\cO^E(\bC).
$$
Hence $\Phi^{-1} g_n-f_n=\Phi^{-1} g_n-\Phi^{-1} e_n=h_n$ is entire, and with any $r>0$
$$
\int_{|\zeta|=r} \Phi^{-1} (\zeta) g_n(\zeta) d\zeta=\int_{|\zeta|=r} f_n (\zeta) d\zeta=(0,\ldots,2\pi n^n i,0,\ldots).
$$
If $\gamma_n$ denotes the $n$th component of $g_n$, it follows that 
$$
\max_{|\zeta|=r}\ {n^{n+1}\over |\zeta|^{n+1}}|\gamma_n (\zeta)|\geq {n^n\over r},\qquad\text{ and }\qquad\max_{|\zeta|=r} \|g_n(\zeta)\|\geq {r^n\over n}.
$$
Choosing $r>1$, we see that $g_n\not\rightarrow 0$ in the compact--open topology, therefore not in the countable cover topology.
This being so for arbitrary lifts $\bg_n$ of $\pi\be_n$, indeed $\pi\be_n\not\rightarrow 0$ in $\Gamma(\bC,\cS)$ as claimed.
\enddemo

\head 7. Completeness\endhead

Whether the $\tau^\delta$ topology on the \v Cech groups $\check H^q(M,\cS)$ is complete is an important question, but one that is unlikely to have a very general answer. In this last section we consider the more modest problem of $\cO^F(M)$, endowed with the $\tau_\delta $ topology,  $F$ a Banach space. That $\cO^F(M)$ is complete is known only when $M$ is a balanced open subset of a Banach (or even Fr\'echet) space, see [D2, Corollary 3.53]. (Dineen formulates the result for $F=\bC$ only, but the proof carries over to any $F$.) The somewhat weaker property of quasi-completeness of $\cO^F(M)$ is known in greater generality. If $M$ is any open subset of a Banach space $X$, then the so--called $\tau_\omega$ topology on $\cO^F(M)$ is complete by [M1, Theorem 5.1]. Further, if $X$ is separable, by [GM, Theorem 2.1] the $\tau_\delta$ and $\tau_\omega$ topologies agree on $\tau_\delta$-bounded subsets of $\cO^F(M)$. (Again, both [M1] and [GM] deal only with $F=\bC$.) It follows that in this case bounded closed subsets of $\big(\cO^F(M),\tau_\delta\big)$ are complete, i.e., $\big(\cO^F(M),\tau_\delta\big)$ is quasi--complete.---I am grateful to Dineen and Mujica for pointing me to [M1, GM].

Although the arguments in [M1, GM] carry over to manifolds as well, for the sake of completeness here we write down the proof of a yet weaker completeness result, which is sufficient for applications we have in mind.

\proclaim{Theorem 7.1} If $M$ is a second countable complex manifold and $F$ is a Banach space, then $\big(\cO^F(M),\tau_\delta\big)$ is sequentially complete.\endproclaim

We shall need the following lemma, whose proof is essentially borrowed from [D1, Proposition 2.4]:

\proclaim{Lemma 7.2} Let $M$ be a complex manifold, $F$ a Banach space, and $f_n\in\cO^F(M)$, $n=1,2,\ldots$, a uniformly bounded sequence. If the $f_n$ form a Cauchy sequence in the $\tau_\delta$ topology, then $f_n$ is locally uniformly convergent.\endproclaim

\demo{Proof} Since this is a local result, we can assume that $M$ is the unit ball in a Banach space $(X,\|\,\|)$, and it will suffice to show that $f_n$ is uniformly convergent on $M'=\{x\in X\: \|x\|<1/2\}$.

Any $f\in \cO^F(M)$ can be expanded in a series of homogeneous polynomials
$$
f=\sum_{k=0}^\infty f^k,\qquad f^k(x)=\int_0^1f(e^{2\pi it}x)e^{-2k\pi it} dt.
$$
If $f$ is bounded, then $\sum f^k(x)$ converges to $f(x)$ uniformly for $\|x\|<r<1$. Indeed, restricting to one dimensional subspaces gives
$$
\big\| f(x)-\sum_{k=0}^p f^k(x)\big\|_F\le\frac{\|x\|^{p+1}}{1-\|x\|}\sup_{\|y\|<1}||f(y)||_F, \qquad ||x||<1.\tag7.1
$$
Furthermore, if $0<r,s\le 1$,
$$
\sup_{||x||<r}||f^k(x)||_F=\big(\frac r s\big)^k\sup_{||x||<s}||f^k(x)||_F\le \big(\frac r s\big)^k \sup_{||x||<s}||f(x)||_F.\tag7.2
$$
Letting $r=1$, it follows that for any $k$ and any countable open cover $\fU$ of $M$ the linear map 
$\cO^F(\fU)\ni f\mapsto f^k\in\cO^F(\{M\})$ is continuous, hence so is the map
$$
\cO^F(M)\ni f\mapsto f^k\in\cO^F(\{M\}).
$$

It follows that for each $k$ the homogeneous components $f^k_n$ of our $f_n$ form a Cauchy sequence in the Banach space 
$\cO^F(\{M\})$, let $\lim_n f_n^k=g^k$. Choosing a positive $A$ so that $\sup_M\|f_n\|_F\le A$ for every $n$,  (7.2)  implies with $r<s=1$
$$
\sup_{||x||<r}\big\|g^k(x)\big\|_F\le Ar^k,\tag7.3
$$
and so $\sum_k g^k=g\in\cO^F(M)$. If $||x||<1/2$ then (7.1), (7.3) imply for any $p$
$$
||f_n(x)-g(x)||_F\le\big\|\sum_{k=0}^p f_n^k(x)-\sum_{k=0}^pg^k(x)\big\|+2^{1-p}A\to 2^{1-p}A,
$$
as $n\to\infty$. Hence $f_n$ indeed converges to $g$, uniformly on $M'$.
\enddemo

\demo{Proof of Theorem 7.1} Let $f_n\in\cO^F(M)$ be a Cauchy sequence in the $\tau_\delta$ topology. The compact--open topology being coarser than $\tau_\delta$, $f_n$ is Cauchy in the compact--open topology as well. In particular, $f_n|K$ are uniformly bounded if $K\subset M$ is compact. This implies that each $x\in M$ has a neighborhood $U$ on which the $f_n$ are uniformly bounded (for otherwise there would be sequences $x_j\in M$, $x_j\to x$, and $n_j\to\infty$ such that $||f_{n_j}(x_j)||_F>j$; but then the $f_n$ would not be uniformly bounded on the compact set $K=\{x,x_1,x_2,\ldots\}$). By Lemma 7.2 the $f_n$ converge locally uniformly to an $f\in\cO^F(M)$.

Thus $M$ is covered by open sets $U$ each on which $f_n\to f$ uniformly. We can select such a countable cover $\fU$; then $f_n\to f$ in 
$\cO^F(\fU)$, and even more in $\big(\cO^F(M),\tau_\delta\big)$. \enddemo

\enddocument
\bye

\bye